\documentclass[reqno, 12pt]{amsart}

\usepackage{amssymb}
\usepackage{latexsym}
\usepackage{epsfig}
\usepackage{epsf}
\usepackage{float}
\usepackage{a4wide}
\usepackage{hyperref}

\newtheorem{theorem}{Theorem}
\newtheorem{lemma}{Lemma}
\newtheorem{corollary}{Corollary}
\newtheorem{proposition}{Proposition}

\newtheorem{remark}{Remark}

\def\bR{{\mathbb R}}  
  
\def\bP{{\mathbb P}}    
\def\bZ{{\mathbb Z}}

\def\bE{{\mathbb E}}

\def\0{{\mathbf 0}}

\def\calP{{\mathcal P}}
\def\calX{{\mathcal X}}
\def\calT{{\mathcal T}}
\def\calD{{\mathcal D}}
\def\calN{{\mathcal N}}

\def\reff#1{(\ref{#1})}
\def\proofof #1{{\noindent \bf Proof of #1.}}
\def\endproof{$\square$ \vskip 2mm}

\begin{document}

\title[On the collision between two PNG droplets]{On the collision between two PNG droplets}
\author{Cristian Coletti and Leandro P. R. Pimentel}
\email{leandro.pimentel@epfl.ch}
\urladdr{http://ima.epfl.ch/$\sim$lpimente/}
\email{cristian@ime.usp.br}
\urladdr{http://www.ime.usp.br/$\sim$cristian}

\begin{abstract}
In this article we study the interface generated by the collision between two cristals growing layer by layer on a one-dimensional substrate through random decomposition of particles. We relate this interface with the notion of $\beta$-path in an equivalent directed polymer model and, by using asymptotics results from Baik and Rains \cite{br2} and some hydrodynamic tools introduced by Cator and Groeenenboon \cite{cg1}, we derive a law of large numbers for such a path and obtain some bounds for its fluctuations.
\end{abstract}

\maketitle

\section{Introduction}
A variety of one dimensional growth models has been proposed to understand the interplay between the geometry of the initial macroscopic profile and the scaling properties of the growing interface \cite{ks}. 
A less well understood phenomen is the interface generated by the collision between two growing materials, named the \emph{competition interface}\footnote{We follow the terminology introduced by Ferrari and Pimentel \cite{fp}}. Since the numerical simulations perfomed by Derrida and Dickman \cite{dd} it is well known that the large space and time behavior of this interface strongly depends on the geometry of the initial profile (see also \cite{sk}). Later, Ferrari, Martin and Pimentel \cite{fmp} considered the competition interface between two clusters in the lattice last-passage percolation set-up and they established a connection between this interface and the so called second-class particle in the totally asymetric exclusion process. This connecion allowed them to perform formal calculations and obtain analitical solutions for the macroscopic description of the competition interface. 

\begin{figure}[hbt]
\begin{center}
\includegraphics[width=0.5\textwidth]{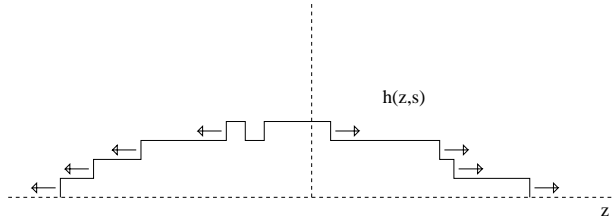}
\end{center}
\caption{Single PNG droplet}\label{1png}
\end{figure}

In this work we do something similar but now in the context of a one dimensional layer by layer growth model  \cite{m}, named the \emph{polynuclear growth} (PNG) model. This model describes a crystal growing layer by layer on a one dimensional substrate through random deposition of particles that nucleate on the existing plateaus of the crystal forming new islands. These islands spread laterally with speed $1$ and adjancent islands of the same level coalesce upon meeting (Figure \ref{1png}).

\begin{figure}[hbt]
\begin{center}
\includegraphics[width=0.5\textwidth]{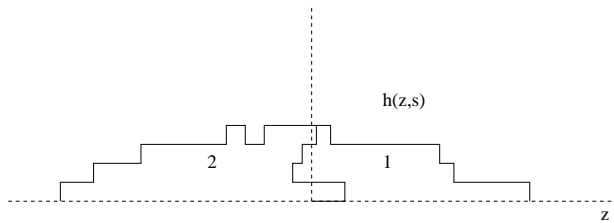}
\end{center}
\caption{Two PNG droplets}\label{2png}
\end{figure}

To consider a multi-type cristal growth model \footnote{For more informations on multi-type growth models we adress to \cite{p} and the references therein.} we assume that the initial substrate is divided into two different types of crystals, say type $1$ if $z<0$ and type $2$ if $z>0$. The dynamics stipulates that if a nucleation occurs on a existing plateau of type $j\in\{1,2\}$ then the new island will be of the same type. When edges of islands having different types meet they stop (Figure \ref{2png}). 

Of course the behavior of this model depends on the geometry of the nucleation events, which can be seen as a point process $\calN$ in $(z,s)$-space-time. We restrict our attention to a particular class of point processes $\calN$ (PNG growth with external sources) and we show a law of large numbers for the competition interface.

The PNG model can be studied in a directed polymer context \cite{ps1} and we show that the space-time path of the competition interface is a particular example of a \emph{$\beta$-path}. The collection of $\beta$-paths form a large class of paths in the directed polymer model and we prove a general theorem that ensures its almost sure  convergence. We also show that the exponent whose value measures the order of the fluctuations of a $\beta$-path about its asymptotic value is at most $2/3$. The proofs are based on the notion of maximal paths and its relation with $\beta$-paths, together with some bounds for the tail of the lenght of the longest directed polymer connecting two distinct points on the space-time plane obtained by Baik and Rains \cite{br2}. 

The PNG and the directed polymer models are intrisictly related to the Hammersley interacting particle system and  another examples of $\beta$-paths arise naturally in this context, named the second-class particle. As a consequence, we obtain a law of large numbers for this special particle which, together with some hydrodynamic ideas introduced by Cator and Groenemboom \cite{cg1}, will also play and important role in studying the asymptotics of $\beta$-paths.

In the next sections we formally introduce all the models considered here and we state the main results.

\subsection{Polynuclear growth with external sources} 
The surface at time $s\geq 0$ is described by an integer-valued function $h(.,s):z\in\bR\to h(z,s)\in\bZ$, named the height profile at time $s$, for which the discontinuity points have upper limits. We consider the initial condition $h(0,z)=0$ for all $z\in\bR$. For each $s>0$ the function $h(.,s)$ has jumps of size one at the discontinuity points, called \emph{up-step} if $h$ increases and \emph{down-step} if $h$ decreases (Figure \ref{1png}). A \emph{nucleation} event at position $(z,s)$ is a creation of a spike, a pair of up- and down-steps, over the previous layer. The up-steps move to the left with unit speed and the down-steps move to the right with unit speed. When an up- and down-steps collide they disappear.

\begin{figure}[hbt]
\begin{center}
\includegraphics[width=0.5\textwidth]{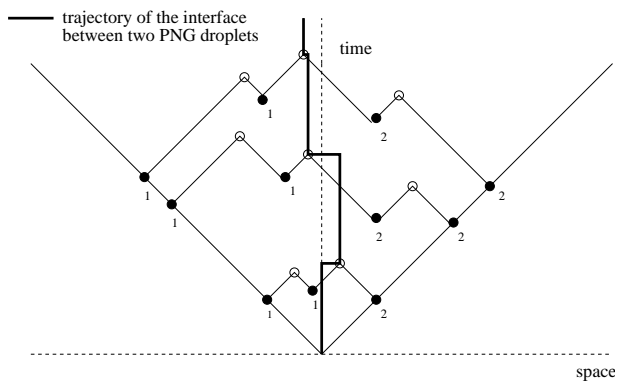}
\end{center}
\caption{The nucleation events. Each point is labeled by its type.}\label{sp-ti}
\end{figure}

The nucleation events form a locally finite point process in space-time. On $\{z=s\}$ we put a Poisson point process $\calD_+$ of intensity $\lambda\geq 0$ while on $\{z=-s\}$ we put a Poisson point process $\calD_-$ of intensity $\rho\geq 0$. On $\{|z|<s\}$ we put a Poisson point process $\calD$ of intensity $1$. We assume that outside $\{|z|\geq s\}$ there is no nucleation event and that all the Poisson point processes involved in this construction are mutually independent. To know the value of $h(z,s)$ one draws the trajectories of the up- and down-steps in the space-time $(z,s)$-plane. When two of these paths meet (as $s$ increases) they stop, which reflects the disappearing of the corresponding up- and down-step. In this way the space-time is divided into regions bounded by piecewise straight lines with slopes equal to $1$ or $-1$ (Figure \ref{sp-ti}). For fixed $z\in\bR^2$ the height $h(z,s)$ is constant in each region.

To introduce the multi-type growth model we assume that the initial layer is divided into two different types of crystals, say $1$ if $z<0$ and $2$ if $z>0$. Consider the  rule which stipulates that if a spike is created over crystal $j$ then it will belong to this material, and that when a down-step of type $1$ collides with an up-step of type $2$ they stop (Figure \ref{2png}). Thus an interface between these two growing crystals is formed. We denote by $\varphi_n$ the position in the $z$-axis of the collision between the different types of crystals at the $n$-th layer and by $\sigma_n$ the time for which this happens, with the convention that $\varphi_0=0$ and $\sigma_0=0$. Define the process $(\varphi(s),\,s\geq 0)$ by setting $\varphi(s)=\varphi_n$ for $s\in[\sigma_n,\sigma_{n+1})$. We prove:
\begin{theorem} \label{png}
Assume that $0\leq\lambda\rho\leq 1$. One has almost surely that
\[
\lim_{s \rightarrow \infty}\frac{\varphi_{s}}{s}= W\,.
\]
where $W\in\left[ \frac{\rho^2-1}{\rho^2 +1},\frac{1-\lambda^2}{1+\lambda^2}\right]$.
\end{theorem}
\begin{remark} \label{remarkpng}
In the stationary regime $\lambda \rho =1$ we obtain a deterministic limiting
value for the inclination of the competition interface. Namely,
$W=\frac{1-\lambda^2}{1+\lambda^2}.$ In Section $2$ we discuss the
correspondence between the limiting value $W$ and the macroscopic behavior of
the height profile.
\end{remark}
\subsection{Directed polymer set-up and $\beta$-paths} 
There is a link between the PNG growth model and a model for directed polymers on Poisson points \cite{ps}. This directed polymer model can be regarded as a last-passage percolation model on $\{x\geq 0\, ,\,t\geq 0\}$ and is defined as follows. Put a Poisson point process $\calP$ of intensity $1$ in the strictly positive quadrant $\{x>0\, ,\,t>0\}$. Independently of $\calP$ we also have mutually independent Poisson point processes, say $\calX$ and $\calT$, on the $x$- and $t$-axis and of intensities $\lambda,\rho\geq 0$, respectively. For $P$ and $Q$ on the plane, define that $P\prec Q$ if both coordinates of $P$ are lower or equal than those of $Q$. For a given realization of the three Poisson point processes, a weakly up/right path, or directed polymer, $(P,P_1,\dots,P_{l},Q)$, starting at $P$ and ending at $Q$, is an oriented and piecewise linear path $\gamma$ connecting $P \prec P_{1} \prec ... \prec P_{l} \prec Q$, where each $P_{j}$ is a Poissonian point for $j=1,\dots,l$. The length $l(\gamma)$ of the  path is the number of Poissonian points used by $\gamma$ and $\Gamma (P,Q)$ denotes the set of all weakly up/right paths from $P$ to $Q$. The maximal length, or the last-passage time, between $P$ and $Q$ is defined by
\begin{equation}
L_m(P,Q) = \max_{\gamma \in \Gamma (P,Q)} l \left( \gamma \right)\,.
\end{equation}
Every $\gamma \in \Gamma (P,Q)$ such that $l(\gamma)=L_{m}(Q)$ is called a maximal path. We aslo consider the growth process $(G_k)_{k\geq 0}$ defined by
\[
G_k=\{Q\,:\,L_m(0,Q)\leq k-1\}
\]
for $k \geq 1$ and for convinience we set $G_0:=\{(0,0)\}$. We denote by $\partial G_k$ the right (hand-side) boundary of $G_k$.
 
Consider the transformation $A\,:\,(x,t)\to (z,s)$ that rotates the $(x,t)$-plane by $45^o$ in the anti-clockwise orientation. If the Poisson point processes involved in the construction of both processes are related by $A$,  $A(\calX)=\calD_-$, $A(\calT)=\calD_+$ and $A(\calP)=\calD$, then the link is apparent. In fact, $h(z,s)$ equals the number of lines crossed by any piecewise linear path from $(0,0)$ to $(z,s)$, with slope between $-1$ and $1$. In particular, one considers the paths which cross them at the nucleation points. These are maximal paths introduced above up to a $45^o$ rotation, and thus it follows that $h(z,s)=L_m(0,(x,t))$.

To see the rule, in this directed polymer model, played by the competition interface we introduce the notion of \emph{$\beta$-paths}. A \emph{$\beta$-point}\footnote{We follow the terminology introduced by Groenemboom \cite{g2} in the Hammersley process context (see Section \ref{Hammersley}).} $Q\in\bR^2_+$ is a concave corner of $\partial G_k$ for some $k\geq 1$. Thus, in the PNG model, $\beta$-points will corresponds to the collisions between up- and down- steps, up to a $45^o$ rotation. We define that $(P_n)_{n\geq 0}$, a sequence of points in $\{x\geq 0\,,\,t\geq 0\}$, is a $\beta$-path if it satisfies: i) $P_n\prec P_{n+1}$; ii) $P_n\in\partial G_n$; iii) $P_n$ is a $\beta$-point. Recalling that $\varphi_n$ denotes the position in the $z$-axis of the collision between the different types of crystals at the $n$-th layer, and that $\sigma_n$ denotes the time in which this happens, one can see that the path $(R_n)_{n\geq 0}$, where $R_{n}:=(\varphi_{n},\sigma_{n})$, is a $\beta$-path up to a $45^o$ rotation (see the trajectory in Figure \ref{sp-ti}). 

For $P=|P|(\cos\theta,\sin\theta)$ and $Q=|Q|(\cos\alpha,\sin\alpha)$, with $\alpha,\theta\in[0,\pi/2]$, let $ang(P,Q)=|\beta-\alpha|$ be the angle in $[0,\pi/2)$ between $P$ and $Q$. We prove: 
\begin{theorem} \label{path}
Assume that $0\leq\lambda\rho\leq 1$. One has almost surely that, if $(P_n)_{n\geq 1}$ is $\beta$-path then
\[ 
\exists\,\lim_{n \rightarrow \infty}\frac{P_{n}}{|P_{n}|}=V=(\cos\theta,\sin\theta)\,, 
\]
where $tan(\theta) \in [\lambda^2,\rho^{-2}]$. 
\end{theorem}

We remark that Theorem \ref{png} follows directly from Theorem \ref{path}. Concerning the fluctuations around its asymptotic angle we have:

\begin{theorem}\label{fluctuation}
Assume that $0\leq\lambda\rho<1$. Then for all $\delta\in(0,1/3)$ there exists a constant $c>0$ such that, almost surely,
\[
ang(P_{m},V) \leq c |P_{m}|^{-\delta}\mbox{ for all large }m\,.
\]
\end{theorem}

We note that Theorem \ref{fluctuation} tell us that, in the regime $0\leq\lambda\rho<1$, for all $\epsilon>0$ the fluctuations of a $\beta$-path $(P_n)_{n\geq 0}$ about its asymptotic value $V|P_n|$ are at most of order $|P_n|^{2/3+\epsilon}$. We do believe that Theorem \ref{path} is almost optimal, i.e. that the correct exponent should be $2/3$.

\subsection{Hammersley process and second-class particles}\label{Hammersley}
Aldous and Diaconis \cite{ad} introduced a continuous time version of the interacting particle process in Hammersley \cite{h} using the following rule. Start with the Poisson point process $\calP$ on $\{x>0\,,\,t>0\}$, of intensity $1$, and move the interval $[0,x]$ vertically through a realization of this point process; if this interval catches a point that is to the right of the points caught before, a new point (or particle) is created in $[0,x]$ at this point; otherwise we shift to this point the previously caught point that is immediately to the right and belongs to $[0,x]$. The number of particles, resulting from this rule, at time  $t$ on the the interval $[0,x]$ is denoted by $N(x,t)$ and the evolving particle process $(N(.,t)\,, t\geq 0)$ is called the \emph{Hammersley process}. In this work we consider an extension of the Hammersley process, as introduced by Groenemboom \cite{g2}, where we also have two others Poisson point processes $\calX$ and $\calT$, of intensities $\lambda$ and $\rho$ and  on the $x$- and $t$-axis, respectively. Points in $\calX$ are called sources while points in $\calT$ are called sinks. Now we have the following rule: start the interacting particle process with a configuration of sources on the $x$-axis, which are subjected to the Hammersley interacting rule in the strictly positive quadrant and which escape through the sinks on the $t$-axis, if such a sink appears to the immediate left of a particle. Now, $N(x,t)$ is the number of particles in $(0,x]\times \{t\}$ plus the number of sinks in $\{0\}\times [0,t]$. When $\rho=1/\lambda$, we have a stationary process \cite{g2}. 

\begin{figure}[hbt]
\begin{center}
\includegraphics[width=0.5\textwidth]{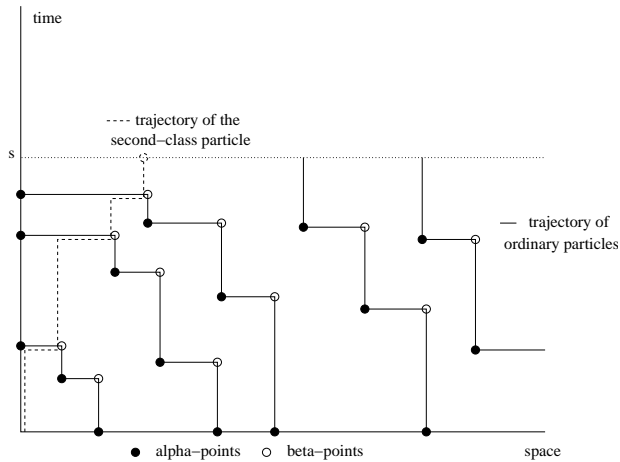}
\end{center}
\caption{Second-class particle}\label{2class}
\end{figure}

Denote by $\Delta_0,\Delta_1,\Delta_2\dots$ the space-time paths of the Hammersley particles with the convention that $\Delta_0=\{(0,0)\}$ and that $\Delta_k$ lies below $\Delta_{k+1}$. Thus, $\partial G_k$ equals to $\Delta_k$ (recall we have constructed both process with the same Poissonian points). Again, if the Poissonian process are related  by $A$, the rotated space-time paths of the up- and down steps correspond to the space-time paths of Hammersley particles. With this picture in mind, one can also see that the $\beta$-points are the left turns of the space-time paths of the particles in the Hammersley process (Figure \ref{2class}). We remark also that, in the stationary regime $\lambda\rho=1$, Cator and Groenemboom \cite{cg1} proved that the $\beta$-points inherit the Poisson property of $\calP$, which allows us to see a duality between $\beta$-paths and maximal paths: a finite $\beta$-path is a maximal path for the time reversal process.

It turns out that another example of a $\beta$-path appears naturally in the Hammersley process: the so called \emph{second-class particles}. A normal second-class particle is a special particle that starts at the origin and jumps to the previous position of the ordinary Hammersley particle that exits through the first sink at the time of the exit, and successively jumps to the previous position of particles directly to the right of it, at times where these particles jump to a position to the left of the second-class particle (Figure \ref{2class}). The position of the second-class particle at time $t$ is denoted by $X_t$. Thus, if $\tau_n$ denotes the time of the $n$-th jump of the second-class particle (with the convention that $\tau_0=0$) then $(Q_n)_{n\geq 0}$, where $Q_n:=(X_{\tau_n},\tau_n)$, is a $\beta$-path. 

\begin{remark}\label{stat} 
In the stationary regime $\lambda\rho =1 $ Cator and Groeneboom \cite{cg1} proved that, almost surely, 
\[
\lim _{t \rightarrow \infty}\frac{X_{t}}{t}=\frac{1}{\lambda^2}\,.
\]
They also showed \cite{cg2} that $X_t-\lambda^{-2}t\sim t^{2/3}$. 
\end{remark}

Here we prove:

\begin{theorem} \label{2par}
Let $(X_t,t \geq 0)$ be the trajectory of a second class particle which is initially at the origin in the Hammersley process with sinks, i.e. $\rho>0$, and such that $\lambda\rho< 1$. Then one has, almost surely,
\[
\lim _{t \rightarrow \infty}\frac{X_{t}}{t}=Z
\]
\noindent where $Z$ is a random variable with the following distribution:
\[
\bP(Z \leq r)=
\left\{
\begin{array}{lll}
0, &r \leq \rho ^{2} ,\\
\frac{\rho ^{-1} - \sqrt{r^{-1}}}{\rho ^{-1}-\lambda}, & \rho ^{2}  < r \leq \lambda ^{-2} ,\\
1, & \lambda ^{-2} < r\,.
\end{array}
\right.
\]
\end{theorem}

The almost sure convergence in the regime $\rho>0$ and $0\leq\lambda\rho<1$ follows directly from Theorem \ref{path}, since the space-time path of a second-class particle can be regard as a $\beta$-path. The description of the limit distribution \footnote{We remark that the limit in distribution of the second-class particle when $\lambda\rho<1$ was also identified by Cator and Dobrynin \cite{cd}.} is obtained in Section \ref{distribution}.

In a further paper \cite{cp} we shall study these models in the regime $\lambda\rho>1$ and we shall prove the almost sure convergence of an arbitrary $\beta$-path to the limit value $\rho/\lambda$. Differently from the regime $\lambda\rho\leq 1$, in this case the fluctuations of $\beta$-paths should be Gaussian. We note that, in the PNG context, this corresponds to the convergence of the competition interface to $(\rho-\lambda)/(\rho+\lambda)$. We also remark that, analogously to the second class particle in the totally asymmetric exclusion process \cite{f1}, the convergence to a deterministic limit value with Gaussian fluctuations is due to the development of a shock in the evolution of the macroscopic profile (hydrodynamic limit).

\subsection*{Overview}
The paper is organized as follows. We begin by studying the stationary regime $\lambda\rho=1$ (Section \ref{Sstat}) and we prove Theorem \ref{path} (in this regime) by using Remark \ref{stat} together with the concept of dual second-class particles. After that we relate the asymptotics for competition interfaces and second-class particles with the respective partial differential equations associated to the macroscopic evolution of the systems. In Section \ref{Sfan}, we start by deriving the convergence in distribution of the second-class particle with coupling ideas of Ferrari and Kpnis \cite{fk} and general hydrodynamics results of Sepp\"al\"ainen \cite{t}. Next we use some results of Baik and Rains \cite{br2} concerning the tail of $L_m$, and the notion of $\delta$-straightness of maximal paths introduced by Newman \cite{n}, to prove the almost sure convergence of $\beta$-paths in the regime $0\leq \lambda\rho <1$ and to obtain the fluctuation upper bound.

\section{Stationary growth and macroscopic description}\label{Sstat}

\subsection{Dual second-class particle} 
\begin{figure}[hbt]
\begin{center}
\includegraphics[width=0.5\textwidth]{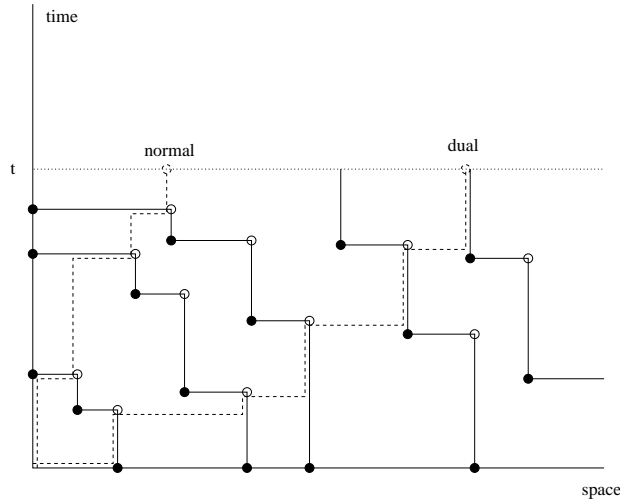}
\end{center}
\caption{Normal and dual second-class particles}\label{b2class}
\end{figure}

The concept of a dual second-class particle was introduced by Cator and Groeneboom \cite{cg1} to prove the convergence of the normal second-class particle in the stationary regime. Recall that to determine the process $t\to L(.,t)$ at point $x$ we shift until time $t$ the interval $[0,x]$ vertically through a realization and we follow the Hammersley interacting rule allowing particles to escape through the sinks. By symmetry, we can also introduce the dual process $x\to L^*(x,.)$ by running the same rule, but now from left to right, i.e. sinks for $L$ becomes sources for $L^*$ and sources for $L$ becomes sinks for $L^*$.  Notice that, in  the stationary regime $\lambda\rho=1$, both processes $L$ and $L^*$ have the same law. We denote  $X^*$ the second-class particle with respect to the dual process $L^*$ and we denote by $X^*_t$ the intersection between the space-time path of the dual second-class particle with $[0,\infty)\times\{t\}$. Trajectories of $X$ and $X^*$ are shown in Figure \ref{b2class}.

\begin{remark}\label{dual} The symmetry of the model and Remark \ref{stat} imply that if $\lambda\rho= 1$ then, almost surely, 
\[
\exists\,\lim _{t \rightarrow \infty}\frac{X^*_{t}}{t}=\frac{1}{\lambda^2}\,.
\]
\end{remark}

An easy but useful observation is that the $\beta$-paths $(Q_n)_{\geq 0}$ and $(Q_n^*)_{\geq 0}$, which correspond to the normal and dual second-class particles are the left- and right-most $\beta$-paths (for all $\lambda,\rho$), respectively:

\begin{lemma}\label{lnd}
Let  $(Q_n)_{\geq 0}$ and $(Q_n^*)_{\geq 0}$ be the $\beta$-paths that correspond to the normal and dual second-class particles respectively, and let $(P_n)_{\geq 0}$ be a $\beta$-path (recall we denote $P=(P(1),P(2))$. Then 
\[
 \frac{Q^*_n(2)}{Q^*_n(1)}\leq \frac{P_n(2)}{P_n(1)}\leq \frac{Q_n(2)}{Q_n(1)}\,.
\]
\end{lemma}

Together with Lemma \ref{lnd}, Remark \ref{dual} implies the convergence of $\beta$-paths in the regime $\lambda\rho=1$.

\subsection{Macroscopic evolution: Hamilton-Jacobi and Burges equations}\label{Macro}
Clearly it is desirable to establish a correspondence between the microscopic structure of the interface generated by the collision between two PNG droplets and its macroscopic behavior. For the PNG droplet, it is known that if $\bar{h}(z,s)$ denotes the \emph{macroscopic height profile} then $\bar{h}$ satisfies the Halmilton-Jacobi equation 
\begin{equation}\label{H-J}    
\partial_s\bar{h}-v(\partial_z\bar{h})=0
\end{equation}
with the \emph{inclination-dependent growth velocity} $v(u)=\sqrt{2+u^2}$ \cite{ps1,t}. For the stationary growth  $\rho\lambda=1$ the solution is $\bar{h}(z,s)= s v(u) + z u$ with $u=(\rho-\lambda)/\sqrt{2}$. Since 
\[
v'(u)=\frac{u}{\sqrt{2 + u^2}}=\frac{\rho-\lambda}{\rho+\lambda}=\frac{1-\lambda^2}{1+\lambda^2}\,
\]
we have that the line $\{z=v'(u) s\}$ is the macroscopic anologue of the competition interface. We also remark that, if one consider the fluctuations of the height profile then the slope $v'(u)$ plays an important rule: the height fluctuations are Gaussian with variance proportional to $t$ except along the line $\{ z=v'(u) s\}$ where they have the KPZ scaling form \cite{ps1,ps}.

In the Hammersley context, we have that if $u(x,t)$ denotes the \emph{macroscopic density profile} then $u$ satisfies the Burgers equation 
\begin{equation}\label{B}
 \partial_t u + \partial_x g(u)\,,
\end{equation}
where $g(u)=1/u$ \cite{t}. The \emph{characteristics} $x(a,t)$ emmanating from $a$  are the solutions to the ordinary differential equation 
\[
x'(t)= g'(u(x,t))\,
\]
with initial condition $x(0)=a$. In the stationary regime $\rho\lambda=1$ the characteristics are given by the lines $x=t\lambda^{-2} + a$, which brings us the macroscopic analogue of the second-class particle, or more generaly, of the $\beta$-paths.

We finish this section by saying a few words concerning the macroscopic evolution when $0\leq \lambda\rho <1$. For the PNG droplet, the solution for the Halmlton-Jacobi equation is of the form $\bar{h}(cs,s)= s f(c)$ where $f$ solves the equation $f(c)-cf'(c)=v(f'(c))$ and the parameter $c$ is related to the local inclination $u=f'(c)$ by $c=-v'(u)$ \cite{ps1, ps}. For instance, when $\lambda=\rho=0$ we have the ellipsoidal shape $f(c)=\sqrt{2(1-c^2)}$. With this information one obtains that the macroscopic height profile has a curved piece between the lines $\{z=s (\rho^2-1)(\rho^2 +1)^{-1}\}$ and $\{z=s (1-\lambda^2)(1+\lambda^2)^{-1}\}$. In the Hammersley context, this corresponds to the development of a \emph{rarefaction front} in the solutions of the Burgers equation, or equivalently, to the existence of infinitelly many characteristics emmanating from the origin. Theorem \ref{path} shows that the macroscopic analogue of a $\beta$-path will be one of these characteristics.

\section{Rarefaction front}\label{Sfan}

\subsection{Convergence in distribution of second class particles}\label{distribution}
The limit law of the second class particle follows from the computation below as well as from Cator and Dobrynin \cite{cd}. Let $\eta _t , t \geq 0$ be the point process obtained by starting with a Poisson Process of intensity $\lambda$ in $(0,\infty)$ at time $0$, and letting it develop according to Hammersley process  on $(0,\infty)$, with Poisson sinks of intensity $\rho$ with $\lambda \rho < 1$ and a Poisson point process of intensity $1$ in the interior of the first quadrant. Furthermore, let $\eta ^{h}_{t} , t \geq 0$ be the process coupled to $\eta _t , t \geq 0$, by using the same points in the first quadrant and on the $t$-axis as used for $\eta$. At time $0$, we consider the same sources on the interval $(h,\infty)$ and on the interval $[0,h]$ we add an independent Poisson process of intensity $\rho ^{-1} - \lambda$. Denote by $\eta _{t}[x,y]$ the number of particles in the interval $[x,y]$ at time $t$ and similarly by $\eta ^{h} _{t}[x,y]$ for the coupled process.

Let
\[
F^{\eta} _{h}(r,t)= \eta _{t}[0,r] - \eta ^{h} _{t}[0,r]\,
\]
and 
\[
F^{\eta}_{h,\epsilon}(r,t)=F^{\eta} _{h}(r \epsilon ^{-1},t\epsilon^{-1})\,. 
\]

Notice that in the absence of extra sources in $[0,h]$ we have $F^{\eta}
_{h,\epsilon}(r,t)=0.$ If there is a unique source in $[0,h]$ coming from the Poisson
point process of intensity $\rho^{-1}-\lambda$ (which happens with
probability $(\rho ^{-1}  - \lambda)h e^{-\rho ^{-1} h}$) we have a discrepancy
which behaves like a second class particle. Denoting it's position at time $t$
by $X_{t,h}$ we get that  $F^{\eta}_{h,\epsilon}(r,t)=-1$ iff $X_{t\epsilon
  ^{-1},h} \leq r\epsilon ^{-1}$. Therefore
\[
\bE\left(F^{\eta}_{h,\epsilon}(r,t)\right)=-\left(\rho ^{-1}  - \lambda \right) h e^{-\rho ^{-1} h} \bP \left(X_{t\epsilon ^{-1},h} \leq r\epsilon ^{-1}\right)+o(h)   
\]
Dividing by $h$ and taking limit when $h$ and $\epsilon$ go to $0$ we get
\[
\lim _{\epsilon \rightarrow 0}\lim _{h \rightarrow 0}\bE\left(\frac{F^{\eta}_{h,\epsilon}(r,t)}{h}\right)=-\left(\rho ^{-1}  - \lambda \right)  \lim _{\epsilon \rightarrow 0}\bP\left(X_{t\epsilon ^{-1}} \leq r\epsilon ^{-1}\right)
\]

On the other hand, by combining the
stationarity of the process $\eta^h$ on $[0,h] \times [0,t]$ with the fact that the number of particles at
time $t$ on $[0,r]$ equals the number of space-time curves crossing the
rectangle $[h,r] \times [0,t]$ plus the number of sources on $[0,h]$, we get  
\begin{align}
\bE\left(F^{\eta}_{h,\epsilon}(r,t)\right)=\nonumber \\
\bE\left(\eta _{t \epsilon ^{-1 }}[r \epsilon ^{-1}-h,r \epsilon ^{-1}]-\eta _{0} ^{h}[0,h]\right)=\nonumber \\
 \bP\left(\eta _{t\epsilon ^{-1}}[r \epsilon ^{-1}-h,r \epsilon ^{-1}]=1,
\eta _{0} ^{h} [0,h]=0\right)-\nonumber \\
\bP\left(\eta _{t\epsilon ^{-1}}[r \epsilon ^{-1}-h,r \epsilon ^{-1}]=0,
\eta _{0} ^{h} [0,h]=1\right) + o(h) \,.\nonumber
\end{align}
Since
\begin{align}
 \bP\left(\eta _{t\epsilon ^{-1}}[r \epsilon ^{-1}-h,r \epsilon ^{-1}]=1,
\eta _{0} ^{h} [0,h]=0\right)-\nonumber \\
\bP\left(\eta _{t\epsilon ^{-1}}[r \epsilon ^{-1}-h,r \epsilon ^{-1}]=0,
\eta _{0} ^{h} [0,h]=1\right) =\nonumber \\
\bP\left(\eta _{t\epsilon ^{-1}}[r \epsilon ^{-1}-h,r \epsilon ^{-1}]=1\right) -
\bP\left(\eta _{0} ^{h} [0,h]=1\right) +
o(h) \,,\nonumber
\end{align}
and 
\[
\bP\left(\eta _{0} ^{h} [0,h]=1\right)=(\rho ^{-1} h) e^{- \rho^{-1} h}\,,
\]
together with the hydrodynamics results of Sepp\"al\"ainen \cite{t}, this yields 
\[
\lim _{\epsilon \rightarrow 0}\lim _{h \rightarrow 0}\bE\left(\frac{F^{\eta}_{h,\epsilon}(r,t)}{h}\right)=u(r,t)-\rho ^{-1}
\]
where $u(r,t)$ is the unique entropic solution of \reff{B}, given by 
\begin{align}
u(x,t)&=\rho ^{-1} \ 1\{\rho ^{-2}x \leq t\}+\sqrt{tx^{-1}} \ 1\{\lambda ^{2} x \leq t < \rho
^{-2} x \}\nonumber \\
&+ \lambda \ 1\{t < \lambda ^{2} x\}\,.\nonumber
\end{align}
Consequently,
\[
\lim _{\epsilon \rightarrow 0} \bP \left(\epsilon X_{t \epsilon ^{-1}} \leq r\right) =
\frac{\rho ^{-1} - u(r,t)}{\rho ^{-1} -\lambda}
\]
which gives the limit law of the second-class particle.

Concerning the limit angle $\theta$ of the competition interface when $0\leq \lambda\rho<1$, at the present moment we do not know how to calculate its distribution but we do know it must be random. In fact, this is true for any $\beta$-path and this is a consequence of the fact that every $\beta$-path is always in between the trajectories of the normal and the dual second-class particles. Indeed, if there exists $a\in [\lambda^2,\rho^{-2}]$ such that with probability one $\tan(\theta)=a$ then with probability one $Z\in[\lambda^2,a]$ and $Z'\in[a,1/\rho^2]$, where $Z$ and $Z'$ denote the limit value of the normal and the dual second-class particles. Although, this leads to a contradiction since we do know that $Z$ and $Z'$ have continuos distributions with support on $[\lambda^2,\rho^{-2}]$.

\subsection{$\delta$-straightness of maximal paths}

In \cite{br2} Baik and Rains used an analytical point of view to study the asymptotics of $L_m$ (see also \cite{bdj} fo the regime $\lambda=\rho=0$). The following is a consequence of their bounds for the tail of $L_m$: if $\lambda_1,\lambda_2\in[0,1)$ then there exists constants $c_j>0$ such that for $k\leq -c_1$ and $t\geq c_2$ we have
\begin{equation}\label{baik1}
\bP(L_{\lambda_1} ^{\lambda_2} (t) - 2t \leq kt^{1/3})\leq c_3 e^{-|k|^3}\,,
\end{equation}
and that for $k\geq c_1$ and $t\geq c_2$ we have
\begin{equation}\label{baik2}
\bP(L_{\lambda_1} ^{\lambda_2} (t) - 2t \geq kt^{1/3})\leq c_4 e^{-c_5
  k^{3/2}}\,,
\end{equation}
where $L_{\lambda_1}^{\lambda_2}(t):=L_m(t,t)$, $\lambda_1$ is the intensity of points in $x$-axis and $\lambda_2$ is the intensity of points in the $t$-axis (see equations (5.2) and (5.14) in \cite{br2}). Notice that for $\lambda\rho<1$ and $\rho^2 < x/t < 1/\lambda^2$,
\begin{align}
\bP(L^{\rho}_{\lambda}(x,t) - 2\sqrt{xt} > k (\sqrt{xt})^{1/3}) &= 
\bP(L^{\lambda_2}_{\lambda_1}(\sqrt{xt}) - 2\sqrt{xt} > k
(\sqrt{xt})^{1/3}) \nonumber
\end{align}
where $\lambda_1=\lambda\sqrt{x/t},\lambda_2=\rho\sqrt{t/x}\in[0,1)$. Therefore, if we denote $\alpha(x,t):=2\sqrt{xt}$, for all $\epsilon>0$ one can find constants $b_j>0$ such that if $\rho^2 +\epsilon< x/t < 1/\lambda^2 -\epsilon$ then 
\begin{equation}\label{baik3}
\bP(|L_m(0,(x,t))-\alpha(x,t)|> k (\sqrt{xt})^{1/3})\leq b_1 e^{-b_2 k^{3/2}}\,.
\end{equation}

Since the work of Kardar, Parisi and Zhang \cite{kpz} it is known there is a strong relation between the fluctuations of $L_m(0,Q)$ and the deviations of a maximal path $\gamma(0,Q)$, connecting $0$ to $Q$, about the line segment $[0,Q]$. It is expected that the scaling relation $\chi=2\xi-1$ holds, if 
\[
|L_m(0,Q)-\alpha(Q)|\sim |Q|^{\chi}\,\mbox{ and }\,\sup_{P\in\gamma(0,Q)}d(P,[0,Q])\sim |Q|^\xi, 
\]
where $d(A,B)$ denotes the euclidean distance between $A$ and $B$. Since for this Poisson last-passage model we do know that $\chi=1/3$, we should have that $\xi=2/3$. Johansson \cite{j} proved that $\xi=2/3$ when $\lambda=\rho=0$ by using a geometric idea developed by Newman \cite{n}, which is based on the curvature properties of the limit shape of the rescaled growth process $n^{-1}G_n$. Newman also introduced the notion of $\delta$-straightness of maximal path as follows. For each Poissonian point $P$, let $R^{out}(P)$ be the set of all Poissonian points $Q$, $P\prec Q$, such that there is a maximal path from $0$ to $Q$ passing through $P$. For $\theta\in(0,\pi/4)$ denote $Co(P,\theta)$ the cone with axis through $P$ and $0$ and with angle $\theta$. Let $\delta>0$. We say that $R^{out}(P)$ is $\delta$-straight if for some constant $c>0$
\[
R^{out}(P)\subseteq Co(P,c|P|^{-\delta})\,.
\]

\begin{proposition}\label{Pstraight}
Assume that $0\leq\lambda\rho<1$. For any $\epsilon>0$ and $\delta\in(0,1/3)$, almost surely, for all but finitely many Poissonian points $P=(P(1),P(2))$ with $P(2)/P(1)\in[\lambda^2+\epsilon,\rho^{-2}-\epsilon]$ one has that $R^{out}(P)$ is $\delta$-straight.
\end{proposition}

For $\lambda=\rho=0$ this is exactly Lemma 2.4 of Wuthrich \cite{w}. To avoid repetitions we give just a sketch of the proof which morally repeat the geometric argument of Newman.

\proofof{Proposition \ref{Pstraight}} Denote by $A_P$ the set of Poisson points $Q$ that satisfies: i) $P\prec Q$; ii) $ang(P,Q)\in[\alpha(P)^{-\delta},2\alpha(P)^{-\delta}]$; iii) $\alpha(Q)\leq 2\alpha(P)$. Notice that $\alpha(aP)=a\alpha(P)$ and so $\alpha(P)$ has the same order of $|P|$. If $|P|$ is sufficiently large then we must have that for all $Q\in A_P$,  $Q(2)/Q(1)\in[\lambda^2+\epsilon/2,\rho^{-2}-\epsilon/2]$. Now, assume there is $Q\in R^{out}(P)\cap A_P$. Then $P$ belongs to some maximal path from $0$ to $Q$ which implies that
\[
L_m(0,Q)=L_m(0,P)+L(P,Q)\,,
\]
and so
\[
\left( \alpha(Q)-L(0,Q)\right)+\left( L(0,P)-\alpha(P)\right) +\left( L(P,Q)-\alpha(Q-P)\right )=
\]
\begin{equation}\label{estraight}
\alpha(Q)-\alpha(P)-\alpha(Q-P)=:\Delta(P,Q)\,.
\end{equation}
By Lemma 2.1 of Wuthrich \cite{w} (which is the desired curvature property for
$\alpha$), for such a $P$ and $Q$,
\[
\Delta(P,Q)\geq |P|^{1-2\delta}\,.
\]
By using \reff{baik3} one can prove that if $\delta\in(0,1/3)$, or equivalently $(1-2\delta)=\chi\in(1/3,1)$, then \reff{estraight} does not occurs for all but finitely many $P$. 

As consequence of the preciding paragraph, one gets that for all but finitely many $P$, if $Q\in R^{out}(P)$ and $\alpha(Q)\leq 2\alpha(P)$ then either
\[
ang(P,Q)\leq \alpha(P)^{-\delta}
\]
or 
\[
ang(P,Q)>2\alpha(P)^{-\delta}\,. 
\]
Since, for sufficiently large $|P|$, to go from $\partial Co(P,\alpha(P)^{-\delta})$ to some point $Q\in Co(P,2\alpha(P)^{-\delta})^c$ a maximal path must pick one Poissonian point $Q'$ with
\[
ang(P,Q')\in[\alpha(P)^{-\delta},2\alpha(P)^{-\delta}]\,,
\]
the second item in the above two possibilities can be delected.

Therefore, for all but finitely many $P$ if $Q\in R^{out}(P)$ and $\alpha(Q)\leq 2\alpha(P)$ then 
\[
ang(P,Q)\leq \alpha(P)^{-\delta}\,. 
\]
Now we claim that this implies $\delta$-straightness. In fact, for every $Q\in Co(P,\epsilon_1)$ the cone $Co(Q,\epsilon_2)$ is contained in the cone $Co(P,\epsilon_1+\epsilon_2)$. By induction, for 
\[
\epsilon_m(P)=\sum_{j=0}^{m-1} (2^j\alpha(P))^{-\delta}\,,
\]
\[
R^{out}(P)\subseteq Co(P,\epsilon_m(P))\bigcup_{\alpha(Q)\geq 2^m \alpha(P)} R^{out}(Q)\,.
\]
By noticing that $\epsilon_m(P)\leq c$, for some constant $c=c(\delta)>0$, one can easily finish this proof. \endproof

As a consequence of the $\delta$-straightness property of maximal paths we have:

\begin{corollary}\label{Cstraight}
Let $a,a'\in(\rho^2,\lambda^{-2})$ with $a<a'$. Almost surely, if $(Q_i)_{i\geq 1}$ and $(Q'_j)_{j\geq 1}$ are two sequences of Poissonian points such that $Q_i\prec Q_{i+1}$,$Q'_j\prec Q'_{j+1}$, $\lim_{i\to\infty}Q_i=\lim_{j\to\infty}Q'_j=\infty$ and 
\[
\limsup\frac{Q'_j(2)}{Q'_j(1)}<1/a'<1/a<\liminf\frac{Q_i(2)}{Q_i(1)}\,
\]
then there are only finitely many $i$ such that, for some $j$, $Q'_j\in R^{out}(Q_i)$. Analogously, there are only finitely many $j$ such that, for some $i$, $Q_i\in R^{out}(Q'_j)$.   
\end{corollary}

\proofof{Corollary \ref{Cstraight}} Divide the positive quadrant into 5 regions as follows: 
\[
C_0:=\{0\leq t\leq \lambda^2 x\}\,,
\]
\[
C_1:=\{0\leq \lambda^2 x\leq t\leq x/a\}\,,
\]
\[
C_2:=\{0\leq x/a\leq t\leq x/a'\}\,, 
\]
\[
C_3:=\{0\leq x/a'\leq t\leq \rho^2\}\,,
\]
and finally,
\[
C_4:=\{0\leq x/\rho^2 \leq t\}\,. 
\]
Pick a $\delta\in(0,1/3)$ and notice that, almost surely, for sufficiently large $|Q|$:
\begin{enumerate}
\item $R^{out}(Q)$ is $\delta$-straight;
\item If $Q\in C_0$ and $Q \prec Q'\in C_3$ then every optimal path from $Q$ to $Q'$ has a Poissonian point in $C_1$, and if $Q\in C_4$ and $Q \prec Q'\in C_1$ then every optimal path from $Q$ to $Q'$ has a Poissonian point in $C_3$; 
\item If $Q\in C_1$ then $Co(Q,c|Q|^{-\delta})\cap (C_3\cup C_4)=\emptyset$ and if $Q\in C_3$ then $Co(Q,c|Q|^{-\delta})\cap (C_1\cup C_0)=\emptyset$.
\end{enumerate}

Now, assume that $Q'_j\in R^{out}(Q_i)$, $|Q'_j|,|Q_i|\geq M$ and $Q'_j(2)/Q'_j(1)<1/a'$. If $Q_i\in C_3$  then, by  (1) and (3), $Q'_j\not\in(C_1\cup C_0)$, which yields to a contradiction. If  $Q_i\in Co_4$, by (2), there exists a $\bar{Q}_i\in C_3$ such that $Q'_j\in R^{out}(\bar{Q}_i)$, and so, by (1) and (3), we also get a contradiction. Since, by assumption, $Q_i\in  C_3\cup  C_4$ for all but finitely many $i$, there are only finitely many $i$ such that, for some $j$, an optimal path from $0$ to $Q'_j$ passes through $Q_i$. The same proof works for the analogue case. \endproof

\subsection{Asymptotics for $\beta$-paths}
\begin{figure}[hbt]
\begin{center}
\epsfxsize 2.0in
\epsfbox{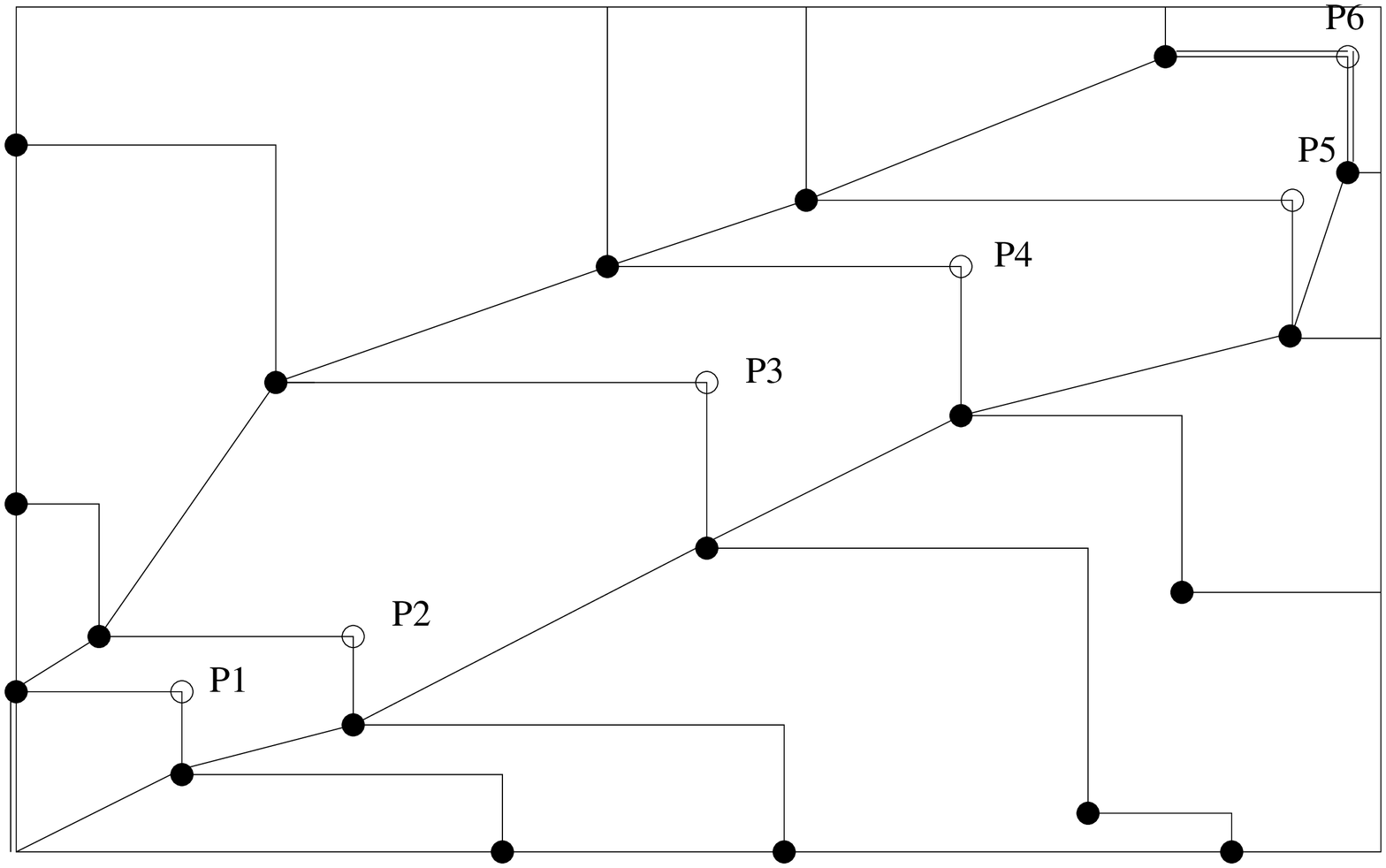}
\caption{\small Two geodesics enclosing a beta-path.}
\label{bpa}
\end{center}
\end{figure}

The idea to control the deviations of a $\beta$-paths, when $0\leq\lambda\rho<1$, is to show that if $(P_n)_{n\geq 0}$ is a $\beta$-path then for all $n\geq 1$ we can construct two maximal paths, both starting from $(0,0)$ and
ending at $P_n$, such that the path $(P_0,\dots,P_n)$ is enclosed by them (see
Figure \ref{bpa}).

\begin{lemma}\label{encl} Almost surely, if $(P_n)_{n\geq 0}$ is a $\beta$-path then for all $n\geq 0$ there exist two maximal paths $\gamma^+_n$ and $\gamma^-_n$ in $\Gamma(0,P_n)$ such that $\gamma^+_n$ is above $(P_0,\dots,P_n)$ and $\gamma^-_n$ is below $(P_0,\dots,P_n)$.
\end{lemma} 

\proofof{Lemma \ref{encl}} Let $G^+_n=(G^+_n(1),P_n(2))$ be the Poissonian point (coming from one of the three Poisson point processes) that first appears to the left (hand-side) of $P_n$ in level $\partial G_n$. Fix $Q\in\bR^2_+$ and let $A_{Q}=\{P \prec Q\}$. Suppose that $G^+_{n},\dots,G^+_{n-k}$ have already been defined for $k< n$. Then set $G^+_{n-(k+1)}$ to be the first Poissonian point in level $\partial G_{n-(k+1)}\cap A_{G^+_{n-k}}$ to the left of $P_{n-(k+1)}$. Notice that if one of the $G^+_{k}$ belongs to the $t$-axis then $G^+_{0},...,G^+_{k-1}$ belong to the $t$-axis as well. By construction, the oriented path $(G^+_0,\dots,G^+_n,P_n)$ is a geodesic (since it picks one point in each level behind $P_n$) which is always above $(P_0,\dots,P_n)$. Similarly, we can construct a geodesic $(G^-_0,\dots,G^-_n,P_n)$ which is below $(P_1,\dots,P_n)$. In this case, we proceed as follows: let $G^-_{n}=(P_n(1),G^-(2)_n)$ be the Poissonian point that first appear to the right of $P_{n}$ in level $\partial G_n$. Suppose that $G_{n}^-,\dots,G_{n-k}^-$ have already been defined for $k<n$. Then set $G^-_{n-(k+1)}$ to be the first Poissonian point in level $\partial G_{n-(k+1)}\cap A_{G^-_{n-k}}$ to the right of $P_{n-(k+1)}$. Notice that if one of the $G^-_{k}$ belongs to the $x$-axis then $G^-_{0},...,G^-_{k-1}$ belong too. By construction, the path $(G^-_0,\dots,G^-_n,P_n)$ is a geodesic which is always below $(P_0,\dots,P_n)$. \endproof

\proofof{Theorem \ref{path} (when $0\leq\lambda\rho<1$)} First we claim that, almost surely, if $(P_n)_{n\geq 1}$ is a $\beta$-path then 
\begin{equation}\label{eenc}
\lambda^2\leq \liminf_{n\to\infty}\frac{P_n(2)}{P_n(1)}\leq\limsup_{n\to\infty}\frac{P_n(2)}{P_n(1)}\leq\frac{1}{\rho^{2}}\,.
\end{equation}   
By Lemma \ref{lnd}, to obtain \reff{eenc1} it suffices to show
\[
\lambda^2\leq \liminf_{n\to\infty}\frac{Q^*_n(2)}{Q^*_n(1)}\leq\limsup_{n\to\infty}\frac{Q_n(2)}{Q_n(1)}\leq\frac{1}{\rho^{2}}\,.
\] 
where $(Q_n)_{\geq 1}$ and $(Q_n^*)_{\geq 1}$ are the $\beta$-paths corresponding to the normal and dual second-class particles, respectively. 

The second inequality follows by coupling the Hammersley process with
parameters $\lambda,\rho$, with the stationary Hammersley process with
parameters $1/\rho,\rho$. Since $1/\rho>\lambda$ (more sources for the
stationary process), $X_t(\lambda,\rho)$ moves to the right faster than
$X_t(1/\rho,\rho),$ i.e. the normal second-class particle for the original
process is always to the right of the normal second-class particle for the
stationary process \cite{cg1}. Together with Remark \ref{dual}, this yields
the second inequality. To show the first inequality, we couple the Hammersley
process, with parameters $\lambda,\rho$, with the stationary process with
parameters $\lambda, 1/\lambda$ (more sinks for the stationary process) and repeat the same argument for the dual second-class particle.

By \reff{eenc}, if $(P_n)_{n\geq 1}$ does not converge then there exist $b<a<a'<b'$ such that
\begin{equation}\label{eenc1}
\lambda^2\leq\liminf_{n\to\infty}\frac{P_n(2)}{P_n(1)}<\frac{1}{b'}<\frac{1}{a'}<\frac{1}{a}<\frac{1}{b}<\limsup_{n\to\infty}\frac{P_n(2)}{P_n(1)}\leq\frac{1}{\rho^{2}}\,.
\end{equation} 
Now let $m<n$ and assume that 
\[
\frac{P_m(2)}{P_m(1)}<\frac{1}{b'} <\frac{1}{b}<\frac{P_n(2)}{P_n(1)}\,. 
\]
Consider the optimal path $\gamma^{-}_n$, giving by Lemma \ref{encl}, which connects $0$ to $P_n$. Since $\gamma^{-}_n$ lies below $(P_0,\dots,P_m,\dots,P_n)$, if $|P_m|$ is sufficiently large then one can find $Q',Q\in\gamma_n^{-}$ (Poissonian points) with $Q\in R^{out}(Q')$ and such that 
\[
\frac{Q'(2)}{Q'(1)}<\frac{1}{a'} < \frac{1}{a} < \frac{Q(2)}{Q(1)}\,. 
\]
Therefore, if \reff{eenc1} occurs then one can construct two sequences of Poissonian points, say  $(Q'_j)_{j\geq 1}$ and $(Q_i)_{i\geq 1}$, with $Q_i\in R^{out}(Q'_j)$ and such that
\[ 
\frac{Q'_j(2)}{Q'_j(1)}<\frac{1}{a'} < \frac{1}{a} < \frac{Q_i(2)}{Q_i(1)}\,
\]
for all $i,j\geq 1$. By Proposition \ref{Pstraight}, this occurs with probability $0$ and thus $(P_n)_{n\geq 1}$ must converge almost surely. \endproof 

\proofof{Theorem \ref{fluctuation}} From Theorem \ref{path}, we have the almost sure convergence of the normal and dual second-class particles and, by a previous calculation (Section \ref{distribution}), their limits have a continuous distribution. Combining this with Lemma \ref{lnd}, one gets that, almost surely, there exists a sufficiently small (random) $\epsilon>0$ such that for all $\beta$-paths $(P_n)_{n\geq 0}$, and $n$ sufficiently large, 
\[
\lambda^2+2\epsilon<\frac{P_n(2)}{P_n(1)}<\rho^{-2}-2\epsilon\,.
\]

Chose a sufficiently large $M$ such that if $|P|\geq M$ and 
\[
\lambda^2+2\epsilon<\frac{P(2)}{P(1)}<\rho^{-2}-2\epsilon\,
\]
then for any $Q\in Co(P,c|P|^{-\delta})$, 
\[
\frac{Q(2)}{Q(1)}\in[\lambda^2+\epsilon,\rho^{-2}-\epsilon]\,. 
\]

Denote by $\theta_m$ the angle in $[0,\pi/2]$ such that
$\tan(\theta_m)=\frac{P_m(2)}{P_m(1)}$ and such  that $|P_m|\geq M$ and that 
\[
\frac{P_n(2)}{P_n(1)}>\tan(\theta_m + 3c|P_m|^{-\delta})\,
\]
for some $n\geq m$. Consider the the maximal path $\gamma_n^-$ giving by Lemma \ref{encl}. Since $\gamma_n^-$ lies below $(P_0,\dots,P_m,\dots,P_n)$, for sufficiently large $M$, there exist $Q,Q'\in \gamma_n^-$ with $Q'\in R^{out}(Q)$ and such that
\[
\frac{Q(2)}{Q(1)}<\tan(\theta_m + c|P_m|^{-\delta})\,\mbox{ and
}\,\frac{Q'(2)}{Q'(1)}>\tan(\theta_m + 2c|P_m|^{-\delta})\,. 
\]
Since $|Q|\sim|P_m|$, this would imply that $R^{out}(Q)$ is not $\delta$-straight which, by Proposition \ref{Pstraight}, occurs with probability $0$.  If 

\[
\frac{P_n(2)}{P_n(1)}<\tan(\theta_m - 3c|P_m|^{-\delta})\,
\]
one can repeat the same argument, but now considering the maximal path $\gamma_n^+$ that lies above $(P_0,\dots,P_m,\dots,P_n)$, to prove that it does not happen with probability $1$. Therefore, for sufficiently large $m$ and for all $n\geq m$,
\[
ang(P_m,P_n)\leq 3c|P_m|^{-\delta}\,.
\]
By sending $n\to\infty$, one gets Theorem \ref{fluctuation}. \endproof


\begin{thebibliography}{60} 

\bibitem{ad} Aldous, D. and Diaconis, P. (1995). Hammersley interacting process and longest increasing subsequences. { \sl Probab. Th. Relat. Fields \bf 103} 199-213.

\bibitem{bdj} Baik, J.,  Deift, P. and Johansson, K. (1999). On the distribution of the lenght of the longest increasing subsequence of random permutations. {\sl J. Amer. Math. Soc. \bf 12} 1119-1178.

\bibitem{br2} Baik, J. and Rains, E. (2000). Limiting distribution for polynuclear growth model with external sources. {\sl J. Stat. Phys. \bf 100} 523-541.

\bibitem{cd} Cator, E. and Dobrynin, S.. Behavior of a second-class particle in Hammersley process (pre-print).

\bibitem{cg1} Cator, E. and Groenemboom, P. (2005). Hammersley process with sources and sinks. {\sl Ann. Probab. \bf 33} 879-903.

\bibitem{cg2} Cator, E. and Groenemboom, P. (2006). Second class particles and cube root asymptotics for Hammersley process. {\sl To appear in Ann. Probab.}.

\bibitem{cp} Coletti, C. and Pimentel, L.P.R. . Shock fluctuations in the Hammersley process (in preparation).

\bibitem{dd} Derrida, B., Dickman, R. (1991). On the interface between two growing Eden clusters {\sl J. Phys. A \bf 24} 191-193.

\bibitem{f1} Ferrari, P.A. (1992). Shock fluctuations in asymmetric simple exclusion. {\sl Probab. Th. Rel. Fields \bf 91} 81-101. 

\bibitem{fk} Ferrari, P.A. and Kipnis, C. (1995). Second class particle in the rarefaction front. {\sl Annales de L'Inst. Henri Poincare \bf 31} 143-154. 

\bibitem{fp} Ferrari, P.A. and Pimentel, L.P.R. (2005). Competition interfaces and second class particles. {\sl Ann. Probab. \bf 33} 1235-1254.

\bibitem{fmp} Ferrari, P.A., Martin, J.B., Pimentel, L.P.R. (2006). Roughening and inclination of competition interfaces. {\sl To appear in Phys. Rev. E}.

\bibitem{g1} Groeneboom, P. (2001). Ulam's problem and the Hammersley process. {\sl Ann. Probab. \bf 29} 683-690.   

\bibitem{g2} Groeneboom, P. (2002). Hydrodynamical methods for analyzing long increasing subsequences .{\sl J. Comp. Appl. Math. \bf 142} 83-105.

\bibitem{h} Hammersley, J.M. (1972). A few seedings of research. In {\sl Proc. 6th Berkeley Symp. Math. Statist. and Probab. \bf 1} 345-394. 

\bibitem{j} Johansoon, K. (2000). Transversal fluctuations for increasing subsequences on the plane. {\sl Probaba. Theory Related Fields \bf 116} 445-456. 

\bibitem{kpz} Kardar, M., Parisi, G. and Zhang, Y.-Z (1986). Dynamic scaling of growing interfaces. {\sl Phys. Rev. Lett. \bf 56}, 889-892.

\bibitem{ks} Krug, J. and Spohn, H. (1992). Kinetic roughening of growing surfaces. In {\sl Solids far from Equilibrium: Growth, Morphology and Defects} (C. Godreche, ed.) Cambridge Univ. Press.  

\bibitem{m} Meakin, P. (1998). Fractals, scaling and growth far from equilibrium, Cambridge Univ. Press.

\bibitem{n} Newman, C. M. (1995). A surface view of first-passage percolation. In {\sl Proc. Intern. Congress of Mathematicians 1994 \bf 2} (S.  D. Chatterji, ed.), Birkhauser, 1017-1023.

\bibitem{p} Pimentel, L.P.R.. A multi-type shape theorem for FPP models (arXiv:math.PR/0411583).

\bibitem{ps1} Pr\"ahofer, M. and Spohn, H. (2000). Statistical self-similarity of one-dimensional growth processes. {\sl Physica A \bf 279}, 342-352.

\bibitem{ps} Pr\"ahofer, M. and Spohn, H. (2005). Universal distributions for growth processes in 1+1 dimensions and random matrices. {\sl Phys. Rev. Lett. \bf 84}, 4882-4885.

\bibitem{sk} Saito Y. and M\"uller-Krumbhaar M. (1995). Critical Phenomena in Morphology Transitions of Growth Models with Competition {\sl Phys. Rev. Lett. \bf 74} 4325-4328.

\bibitem{t} Sep\"al\"ainen, T. (2002). Diffusive fluctuations for one-dimensional
  totally asymmetric interacting random dynamics {\sl Comm. Math. Phys. \bf
  229}  141-182.

\bibitem{w} W\"{u}thrich, M. V. (2000). Asymptotic behavior of semi-infinite geodesics for maximal increasing subsequences in the plane. In {\sl In and Out of Equilibrium} (V. Sidoravicius ed.), Birkhauser, 205-226.

\end{thebibliography}
\end{document}